\documentclass[12pt]{article}
\usepackage{mathptmx}
\newtheorem{theorem}{Theorem}
\newtheorem{lemma}{Lemma}
\newcommand{\ee}{\mathrm{e}}
\begin{document}
\title{Asymmetric Latin squares, Steiner triple systems, and edge-parallelisms}
\author{Peter J. Cameron\footnote{My address when I wrote this paper was
Merton College, Oxford OX1~4JD, UK. My current address is School of Mathematics
and Statistics, University of St Andrews, North Haugh, St Andrews, Fife
KY16~9SS, UK.}}
\date{}
\maketitle
\begin{quote}
This article, showing that almost all objects in the title are asymmetric,
is re-typed from a manuscript I wrote somewhere around 1980 (after the 
papers of Bang and Friedland on the permanent conjecture but before those of
Egorychev and Falikman). I am not sure of the exact date. The manuscript had
been lost, but surfaced among my papers recently.

I am grateful to Laci Babai and Ian Wanless who have encouraged me to make
this document public, and to Ian for spotting a couple of typos. In the
section on Latin squares, Ian objects to my use of the term ``cell''; this
might be more reasonably called a ``triple'' (since it specifies a row, column
and symbol), but I have decided to keep the terminology I originally used.

The result for Latin squares is in
\begin{itemize}
\item[] B. D. McKay and I. M. Wanless, On the number of Latin squares,
\textit{Annals of Combinatorics} \textbf{9} (2005), 335--344 (arXiv 0909.2101)
\end{itemize}
while the result for Steiner triple systems is in
\begin{itemize}
\item[] L. Babai, Almost all Steiner triple systems are asymmetric,
\textit{Annals of Discrete Mathematics} \textbf{7} (1980), 37--39.
\end{itemize}
\end{quote}

\section{Introduction}

Recently, Bang~\cite{1} and Friedland~\cite{3} have shown that the permanent
of a doubly stochastic matrix of order~$n$ is at least $\ee^{-n}$. This
result substantially improves known lower bounds for the numbers of
combinatorial structures of the types mentioned in the title. (It is already
documented in the literature~\cite{6,8,2} that such improvement would follow
from the truth of the van der Waerden permanent conjecture; the result of
Bang and Friedland is close enough to the conjecture to have the same effect.)
In this paper, I give a possibly less well-known consequence of the result
on permanents.

\begin{theorem}
Almost all Latin squares, Steiner triple systems, or edge-parallelisms of
complete graphs have no non-trivial automorphisms; that is, the proportion
of such objects of an admissible order $n$ admitting non-trivial 
automorphisms tends to zero as $n\to\infty$.
\end{theorem}

Here, as is well-known, $n$ is admissible for Steiner triple systems if and
only if $n\equiv1$ or $3$ (mod~$6$), and $n$ is admissible for
edge-parallelisms if and only if $n\equiv0$ (mod~$2$). All integers are
admissible orders of Latin squares. The paper concludes with the observation
that a similar result holds for strongly regular graphs with least eigenvalue
$-3$ or greater.

I am grateful to J.~H.~van~Lint for helpful discussions on permanents.

\section{Latin squares}

Given an $n\times(n-k)$ Latin rectangle, the number of ways of choosing an
$(n-k+1)^{\mathrm{st}}$ row is the permanent of a $(0,1)$ matrix of order~$n$
with row and column sums $k$ (see Ryser~\cite{6}), and hence is at least
$(k/\ee)^n$ (by~\cite{1,3}). So the number of Latin squares of order $n$ is
at least $\displaystyle{\prod_{k=1}^n(k/\ee)^n=(n!)6/\ee^{n^2}}$. This
number is greater than $n^{(1-\epsilon)n^2}$ for $n\ge n_0(\epsilon)$.

We take the most general definition of an automorphism of a Latin square $S$,
as a permutation on the $3n$ symbols indexing the rows, columns and entries
(say $\{r_1,\ldots,r_n,c_1,\ldots,c_n,e_1,\ldots,e_n\}$) preserving the
obvious partition into three sets $R,C,E$ of size $n$ and also the set of
triples $(r_i,c_j,e_k)$ for which the $(i,j)$ entry of $S$ is $k$. (We call
such triples \emph{cells}.) If an automorphism fixes elements in at least
two of $R,C,E$, then its fixed elements form a subsquare of $S$. Note that
the order of a subsquare is at most $\frac{1}{2}n$.

Now let $g$ be one of the $6(n!)^3$ permutations of $R\cup C\cup E$ fixing
the partition. How many Latin squares admit $g$ as an automorphism? If $g$
doesn't fix the three sets $R,C,E$, then it fixes at most $n$ cells of any
such square (for any fixed cell on $r_i$ must also be on $c_j$, if
$g(r_i)=c_j$, and $r_i$ and $c_j$ determine a unique cell; similar arguments
in the other cases). If $g$ is not the identity but fixes the three sets then,
as remarked earlier, it fixes at most $\frac{1}{4}n^2$ cells. For $n\ge4$,
we have $n\le\frac{1}{4}n^2$.

Let $r$ be the number of fixed cells (determined by their rows and columns).
We may choose their entries in at most $n^r$ ways. Any choice of entry for
a non-fixed cell determines all the cells in its orbit under $g$; so there
are at most $n^{\frac{1}{2}(n^2-r)}$ of these. So the number of fixed
squares is at most $n^{\frac{1}{2}(n^2+r)}\le n^{5n^2/8}$.

Hence the number of Latin squares admitting non-trivial automorphisms is at
most $6(n!)^3n^{5n^2/8}=o((n!)^n/\ee^{n^2})$.

\section{Steiner triple systems}

The number of Steiner triple systems of admissible order $n$ is at least
$n^{(1-\epsilon)n^2/6}$ for sufficiently large $n$ (combining Wilson's
results~\cite{8} with those of Bang and Friedland).

Let $g$ be a non-identity automorphism of a Steiner triple system $S$ of
order $n$, and suppose $g$ fixes $m$ points. The fixed points carry a
subsystem of $S$, so $m\le\frac{1}{2}(n-1)$. This subsystem contains
$m(m-1)/6$ fixed blocks. Any other point lies in at most one fixed block,
so at most $\frac{1}{2}(n-m)$ further blocks are fixed. The total number
of fixed blocks is thus at most $(n^2+2n-9)/24$, and the number $r$ of
block-orbits satisfies
\begin{eqnarray*}
r &\le& (n^2+2n-9)/24+\frac{1}{2}(n(n-1)/6-(n^2+2n-9)/24) \\
  &<& 5n^2/48.
\end{eqnarray*}

Now take a permutation $g$ on the set of points. Choose triples for the blocks
of a Steiner triple system admitting $g$ in such a way that, when any new
block is chosen, its entire orbit under $g$ is included. The number of such
sequences of blocks is at most ${n\choose3}^r<(n^3/6)^r$; so the number of
Steiner triple systems is at most $(n^3\ee/6r)^r$.

Now $(a/x)^x$ is an increasing function of $x$ for $x<a\ee$; so, since
$r\le 5n^2/48$, we have that $(n^3\ee/6r)^r\le(8n\ee/5)^{5n^2/48}$. Hence
the number of Steiner triple systems admitting non-trivial automorphisms
is at most $n!(8n\ee/5)^{5n^2/48}=o(n^{(1-\epsilon)n^2/6})$.

\section{Edge-parallelisms}

The structures considered here are sometimes referred to as $1$-factorisations
or minimal edge-colourings of complete graphs; they are partitions of the
$2$-subsets of an $n$-set $X$ into ``parallel classes'', each of which 
partitions $X$. For a general reference, see~\cite[Chapter 4]{2}. It follows
from~\cite{2} together with the result of Bang and Friedland that, if $n$ is
admissible (that is, even), the number of edge-parallelisms of order $n$ is
at least $n^{(1-\epsilon)n^2/2}$ for $n\ge n_0(\epsilon)$.

We need the fact that the number of $1$-factors of a $k$-valent graph on $n$
vertices is at most $k^{\frac{1}{2}n}$ (see~\cite[p. 64]{2}).

\begin{lemma}
Let $\Gamma$ be a $k$-valent graph on $n$ vertices, $g$ an automorphism of
$\Gamma$ with no fixed vertices. Then the number of $1$-factors of $\Gamma$
fixed by $g$ is at most $(8\ee k)^{\frac{1}{4}n}$.
\end{lemma}

\paragraph{Proof} Count fixed $1$-factors containing $r$ edges fixed by $g$.
The fixed edges are $2$-cycles of $g$, so there are at most $\frac{1}{2}n
\choose r$ choices for these. Suppose the non-fixed edges lie in $m$ orbits
under $g$. Choosing these in order, such that each new edge chosen is followed
by its orbit, we have at most $((\frac{1}{2}n-r)k)^m$ choices; hence at most
$((\frac{1}{2}n-r)k)^m/m! < ((\frac{1}{2}n-r)k\ee/m)^m$ choices up to
permutations of the orbits. As in the last section, this number is greatest
when $m$ has its largest possible value $\frac{1}{2}(\frac{1}{2}n-r)$, and so
it is smaller than $(2\ee k)^{\frac{1}{2}(\frac{1}{2}n-r)}$. Now the total
number of $1$-factors is less than
\[\sum_{r=0}^{\frac{1}{2}n}{\frac{1}{2}n\choose r}(2\ee k)^{\frac{1}{2}(\frac{1}{2}n-r)} \le 2^{\frac{1}{2}n}(2\ee k)^{\frac{1}{4}n} = (8\ee k)^{\frac{1}{4}n}.\]

\medskip

Now we turn to the proof of the theorem. Suppose $g$ is a permutation of an
$n$-set; we want to count edge-parallelisms fixed by $g$. If $g$ fixes $r$
points, with $r>0$, then its fixed points carry a subsystem, whence
$r\le\frac{1}{2}n$ (\cite[p. 25]{2}), and it fixes $r-1$ parallel classes
($1$-factors). So the number of orbits of $g$ on parallel classes satisfies
$m\le r+\frac{1}{2}(n-r)\le\frac{3}{4}n$. There are at most $n^{\frac{1}{2}n}$
$1$-factors altogether, and so at most $n^{3n^2/8}$ fixed edge-parallelisms.

Now suppose that $g$ fixes no points; count fixed edge-parallelisms with $s$
fixed parallel classes. By the lemma, the fixed parallel classes can be chosen
in at most $(8\ee n)^{\frac{1}{4}ns}$ ways. If the remaining classes fall into
$m$ orbits, then $m\le\frac{1}{2}(n-s)$, and as before there are at most
$n^{\frac{1}{4}n(n-s)}$ choices for these. Multiplying, and summing over $s$,
we obtain at most $n(8\ee n)^{\frac{1}{4}n^2}$ fixed edge-parallelisms. This
number is smaller than $n^{3n^2/8}$ for sufficiently large $n$.

Thus the number of edge-parallelisms admitting non-trivial automorphisms is
at most $n!\,n^{3n^2/8}=o(n^{(1-\epsilon)\frac{1}{2}n^2})$.

\section{Strongly regular graphs}

Ray-Chaudhuri~\cite{5} and Neumaier~\cite{4} have shown that all but finitely
many strongly regular graphs with least eigenvalue $-3$ are of one of the
following types:
\begin{itemize}\itemsep0pt
\item[(i)] complete multipartite with block size $3$;
\item[(ii)] a Latin square graph (whose vertices are the cells of a Latin
square, two vertices adjacent if the cells agree in row, column or entry);
\item[(iii)] a Steiner graph (whose vertices are the blocks of a Steiner
triple system, two vertices adjacent if the blocks intersect in a point).
\end{itemize}

For all but finitely many graphs of the second and third type, every
graph-automorphism is induced by an automorphism of the Latin square or
Steiner triple system. Moreover, all but finitely many strongly regular
graphs with least eigenvalue greater than $-3$ are complete multipartite
with block size $2$, or square lattice or triangular graphs (Seidel~\cite{7}).

It follows that, of strongly regular graphs with least eigenvalue $-3$ or
greater on at most $n$ vertices, the proportion admitting non-trivial
automorphisms tends to zero as $n\to\infty$.

It would be interesting to know whether the same assertion holds without
the restriction on the least eigenvalue.


\begin{thebibliography}{8}

\bibitem{1}
T. Bang, On matrix functions giving a good approximation to the van der
Waerden permanent conjecture, preprint no.~30, Copehnagen University, 1979.

\bibitem{2}
P. J. Cameron, ``Parallelisms of Complete Designs'', London Math. Soc.
Lecture Notes \textbf{23}, Cambridge Univ. Pr., Cambridge, 1976.

\bibitem{3}
S. Friedland, A lower bound for the permanent of a doubly stochastic matrix,
\textit{Ann. Math.} \textbf{110} (1979), 167--176.

\bibitem{4}
A. Neumaier, Strongly regular graphs with least eigenvalue $-m$, to appear.

\bibitem{5}
D. K. Ray-Chaudhuri, Uniqueness of association schemes, \textit{Proc. Int.
Colloq. Teorie Combinatorie}, 465--479, Accad. Naz. Lincei, Roma, 1977.

\bibitem{6}
H. J. Ryser, Permanents and systems of distinct representatives,
``Combinatorial Mathematics and its Applications'', 55--68, Univ. North
Carolina Pr., Chapel Hill, 1969.

\bibitem{7}
J. J. Seidel, Graphs and two-graphs, \textit{Proc. Fifth Southeastern Conf.
Combinatorics, Graph Theory, Computing}, 125--143, Congressus Numerantium X,
Utilitas Math., Winnipeg, 1974.

\bibitem{8}
R. M. Wilson, Nonisomorphic Steiner triple systems, \textit{MAth. Z.}
\textbf{135} (1974), 303--313.

\end{thebibliography}
\end{document}